# Group Theory in School Mathematics? Teaching Permutation Cycles Through the 15-Puzzle

Bence Torma[1], Tamás Waldhauser[2]

[1] Doctoral School of Mathematics, University of Szeged, Hungary,

[2] Bolyai Institute, University of Szeged, Hungary,

tormabence97@gmail.com

*Abstract: Permutation cycles are generally associated with undergraduate abstract algebra. This exploratory study examined whether students in Grades 5–11 could construct and use cycle representations in the context of the 15-puzzle. After a 45-minute teacher-guided lesson moving from puzzle manipulation to arrow diagrams and cycle notation, 313 students analyzed one of two new configurations—one solvable and one unsolvable—and used a supplied rule to classify it. Of these students, 78.3% constructed a correct cycle representation, and 67.7% both constructed the representation correctly and reached the correct classification. Cycle-construction accuracy was similar for the two configurations, but classification was less often correct for the unsolvable configuration. The findings concern immediate, supported performance rather than full understanding of permutation cycles or group theory. Nevertheless, they show that many students could use cycle representations after brief instruction and that constructing the representation and using it to reach a conclusion were separate demands. One reason to teach permutations is that, as finite, discrete, non-formulaic functions, they can extend students' experience of functions beyond familiar formulae and continuous graphs.*

## INTRODUCTION

When permutations appear in school mathematics, they are usually treated as arrangements to be counted. Their algebraic structure—and in particular their decomposition into cycles—is generally associated with undergraduate abstract algebra. This raises a deliberately provocative question: can an idea drawn from group theory be taught meaningfully much earlier?

This exploratory classroom study examined whether students in Grades 5–11 could construct and use cycle representations in the context of the 15-puzzle. After a single 45-minute teacher-guided lesson, 245 of the 313 participants (78.3%) produced a correct cycle representation of a new puzzle configuration, and 212 (67.7%) both produced the correct representation and reached the correct solvability classification. These results do not show that the students had acquired a full structural understanding of permutation cycles, still less of group theory. They do show, however, that many of them could construct and use the representation successfully after brief, guided instruction.

Feasibility alone, however, is not a reason to add material to the school curriculum. Why teach permutation cycles at all? Our broader argument is that permutations—and discrete functions more generally—can extend students' experience of functions beyond the formulae, tables, and continuous graphs that dominate school mathematics (Tall & Vinner, 1981; Vinner & Dreyfus, 1989). As bijections of finite sets, permutations can be displayed in full and need not be described by an algebraic rule. Their cycle decompositions also make structural features visible that are difficult to see directly in an arrangement or table. The present study does not test whether a single lesson changes students' concept of function, nor does it compare the puzzle lesson with another teaching approach. It addresses a necessary prior question: whether students can construct and use the representation on which such instruction would depend.

The 15-puzzle offers a natural setting for this question. Its immediate goal is easy to understand: numbered tiles are slid into the empty square until the prescribed order is restored. Students can begin by experimenting with the puzzle itself. The familiar 15–14 configuration then changes the character of the activity. Swapping only tiles 14 and 15 produces a configuration that looks almost solved but cannot be reached from the solved configuration by legal moves. Understanding why requires more than continued trial and error; students need a property of the configuration that is preserved by every move.

In the lesson examined here, the move from experimenting with the puzzle to seeking a mathematical explanation of solvability was supported by a sequence of representations: the physical puzzle, arrow diagrams, and finally cycle notation. The same configuration was viewed in increasingly structural ways, while the practical question of solvability remained unchanged. The literature reviewed below focuses on four ideas directly relevant to the design of the lesson: rich tasks, coordination of representations, structural thinking, and permutations as discrete functions.

The analysis distinguished two kinds of success: constructing a correct cycle representation and using the explicitly taught rule to decide whether the configuration was solvable. These are not equivalent achievements, so they were analyzed separately.

The study addresses the following research questions:

1. How did students in Grades 5–11 respond to the lesson in terms of perceived interest, surprise at its playful character, prior expectations of substantial mathematical background and self-reported understanding?
2. To what extent were students able to construct a correct cycle representation of a puzzle configuration and reach the correct conclusion about solvability?
3. Did performance differ between the solvable and unsolvable task versions?

The remainder of the paper sets out the theoretical and mathematical background of the lesson, describes the classroom study, reports the results, and considers their pedagogical implications and limitations.

## LITERATURE REVIEW

### Rich Tasks, Mathematical Puzzles and Exploratory Learning

The 15-puzzle has several features of a rich mathematical task. It can be entered with little formal knowledge, yet it opens the way to substantially deeper questions. Instead of applying a predetermined algorithm, students can experiment, notice patterns, make conjectures, and look for mathematical structure (Griffiths, 2010). In this sense, the task has a low threshold but a high ceiling: direct play with the puzzle can lead naturally to permutations, cycle decompositions, and invariants (Chamberlin, 2019).

Puzzles are useful for this kind of work because the mathematical question grows out of an intelligible goal. Research on playful and problem-based settings has shown how such contexts can support exploration of discrete ideas and help students develop meaning for combinatorial concepts (Quinn & Wiest, 1998; Simonson & Holm, 2003). In the 15-puzzle, cycle notation is not introduced as an isolated symbolic technique: it appears because students need a way to decide whether a configuration is solvable, thereby embedding the procedure in a mathematical problem with a clear purpose (Foster, 2013).

Evidence from recent classroom work also suggests that non-digital mathematical games can do more than provide practice; they can contribute to confidence, strategic awareness, and mathematical understanding (Fuentes, 2024).

Richness in the task itself is not enough; what teachers do with the task matters as well. Griffiths (2009) stresses the importance of anticipating likely student approaches, offering support at appropriate points, and attending to difficulties that emerge during the lesson. That was also the role of guidance in the present lessons. Students began by experimenting, but the teacher then helped them connect the visible tile arrangement to arrows, cycles, and finally the solvability criterion.

More broadly, problem-based classrooms can give students opportunities to build mathematical knowledge through exploration, contextualized problem solving, and reflection rather than through procedural work alone (Martín-Cudero et al., 2026).

### Representations and the Transition to Structural Thinking

The 15-puzzle is educationally interesting not only because it is playful, but because the same mathematical situation can be represented in several ways. At first, students see and manipulate an arrangement of tiles. To reason about solvability, however, attention has to shift from individual moves to the relationship between current and target positions. This creates a reason to move from the concrete board to more abstract representations.

Duval's framework treats conversion between semiotic registers as an important part of mathematical understanding, not merely the use of several representations side by side (Duval, 2006). Work on multiple representations makes a similar point: their usefulness depends on what each representation contributes and on whether learners can coordinate them (Ainsworth, 2006; González-Martín et al., 2011). In this lesson, the board shows the configuration itself, arrows make the mapping between positions explicit, and cycle notation records that mapping compactly. None of these forms is useful simply by being present; students have to see how they refer to the same mathematical object. This coordination can be difficult, and success with one representation does not necessarily imply success in interpreting another (Pathak, 2025).

Sfard's distinction between operational and structural conceptions provides a useful framework for this progression (Sfard, 1991). Operational reasoning may focus on a sequence of possible moves, whereas a structural perspective treats the whole configuration as one permutation with a global organization. The progression from board to arrows to cycles was intended to support this change of viewpoint. For this reason, the analysis distinguishes the ability to produce the cycle representation from the ability to use the structural information it contains.

Related classroom research in abstract algebra has likewise found that guided exploration and intermediate representations, supported by carefully designed worksheets, can help students move from concrete examples towards formal algebraic ideas (Pokhrel et al., 2026).

### Permutations as Discrete Functions

Permutations also offer an unusual but mathematically straightforward example of a function. A permutation is a bijection from a finite set to itself. Unlike the functions that dominate much school mathematics, it need not be given by a formula and its graph need not be viewed as a continuous curve. Because the domain is finite, the whole mapping can be displayed at once.

The use of permutations as examples of functions can also be viewed through Dienes's principles of mathematical and perceptual variability (Dienes, 1960). Mathematical variability means varying non-essential features while retaining what defines the concept. In the present context, permutations supply finite, discrete, non-formulaic examples of functions. Perceptual variability concerns meeting the same mathematical object in different forms. In the 15-puzzle, the same permutation can be seen as a tile arrangement, a mapping between positions, an arrow diagram, a table, or a cycle decomposition.

This is important because students' concept images of functions can become narrower than the formal definition when most examples involve familiar formulae and real-valued graphs. Students may, for example, expect every function to have an algebraic rule or a continuous-looking graph (Tall & Vinner, 1981; Vinner & Dreyfus, 1989). Permutations keep the defining input-output correspondence while dropping those non-essential features.

The 15-puzzle also gives the mapping a concrete meaning: it records where the elements of the target configuration are currently located. Cycle notation then exposes the structure of that mapping. We do not claim here that one lesson changes students' concept of function. Instead, we examine a prerequisite for any such longer-term claim: whether school students can first learn to construct and use the cycle representation in this setting.

## THE 15-PUZZLE

Permutation puzzles require a set of elements—for example, tiles, numbers, or colored pieces—to be transformed from one arrangement to another under specified legal moves. The 15-puzzle and Rubik's Cube are familiar examples.

The 15-puzzle has fifteen numbered tiles in a 4 × 4 frame, leaving one square empty. In the solved configuration (Figure 1), the tiles appear in numerical order from the top-left corner, row by row, with the bottom-right square empty. A legal move consists of sliding a tile adjacent to the empty square into that square. The goal is to reach the solved configuration using only such moves. In early versions of the puzzle, the tiles could be removed from the frame and replaced in an arbitrary order before play began. Such an arrangement was not necessarily solvable and could therefore produce configurations such as the famous 15–14 configuration (Slocum & Sonneveld, 2006). In most modern versions, by contrast, the tiles remain in the frame, so starting configurations are obtained by sliding the tiles from a solvable configuration and are therefore themselves solvable.

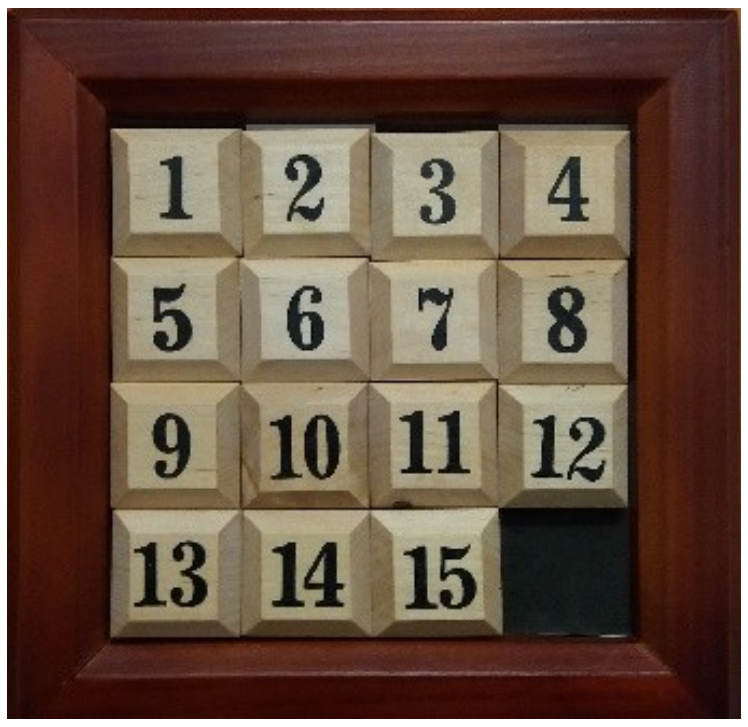

*Figure 1: The solved configuration*

*Note. All photographs and diagrams in this article were created by the authors. The 15-puzzle shown in the photographs is the puzzle supplied with the book The Famous 15 Puzzle (Slocum, 2009).*

### The Puzzle and the Solvability Problem

The puzzle was first marketed in Boston in 1879 and quickly became a craze in the United States and Europe. A particularly important feature of the puzzle is that, for a *solvable* configuration, successfully placing the first thirteen numbered tiles is already sufficient to guarantee completion: once tiles 1–13 have been brought to their correct positions, tiles 14 and 15 must necessarily occur in the correct order. Thus, a solver who had learned how to arrange the first thirteen tiles could reasonably believe that a general solution method had been found.

This made the famous 15–14 configuration especially deceptive. In this configuration, all tiles are correctly placed except that 14 and 15 are interchanged. A solver starting from an *unsolvable* configuration could still arrange the first thirteen tiles correctly, only to arrive at the final 15–14 arrangement. Since the same procedure would always lead to successful completion from a solvable configuration, it was natural to suspect that some mistake had been made rather than to accept that the configuration itself was unsolvable. This helps explain why the 15–14 configuration attracted so much attention and why challenges and monetary prizes were offered for solving it. The puzzle later became closely associated with Sam Loyd, although it predates his involvement (Slocum & Sonneveld, 2006).

For teaching, the 15–14 configuration is interesting for another reason as well. It raises a natural mathematical question: why does a procedure that appears to solve almost the entire puzzle sometimes become blocked only at the very end, and how can we know in advance that no sequence of legal moves can succeed? The approach used in the lesson answers this question by combining checkerboard coloring with the cycle structure of the corresponding permutation to obtain an invariant that gives a condition that is both necessary and sufficient for solvability.

### Representing a Configuration by Cycles

The idea was introduced in the same informal way in all 16 lessons. For this purpose, the empty square was treated as a fictitious sixteenth tile, labeled $E$. Students were asked to think of an arrow from each tile, including $E$, to the tile occupying its target position in the solved puzzle (Figure 2a). Starting from any tile and following these arrows eventually leads back to the starting tile, forming a cycle. The same procedure is then repeated from a tile not yet included in a cycle, until all sixteen tiles, including $E$, have been used (Figure 2b).

For completeness, we give a precise mathematical description of the same construction. Let $X = \{1, 2, \dots, 15, E\}$, where 1, …, 15 are the tile labels and $E$ represents the empty square. For a given puzzle configuration and $x \in X$, define $p(x)$ to be the element of $X$ occupying the target position of $x$. Then $p$ is a permutation of $X$. The arrow $x \mapsto p(x)$ therefore goes from each element to the element occupying its target position. If we start from any $x \in X$ and follow the arrows, we eventually return to $x$ and obtain a cycle. Repeating this procedure for elements not yet included gives the decomposition of $p$ into disjoint cycles. An element in its target position is a fixed point. Throughout this paper, the total number of cycles includes all fixed points, each of which is counted as a one-element cycle.

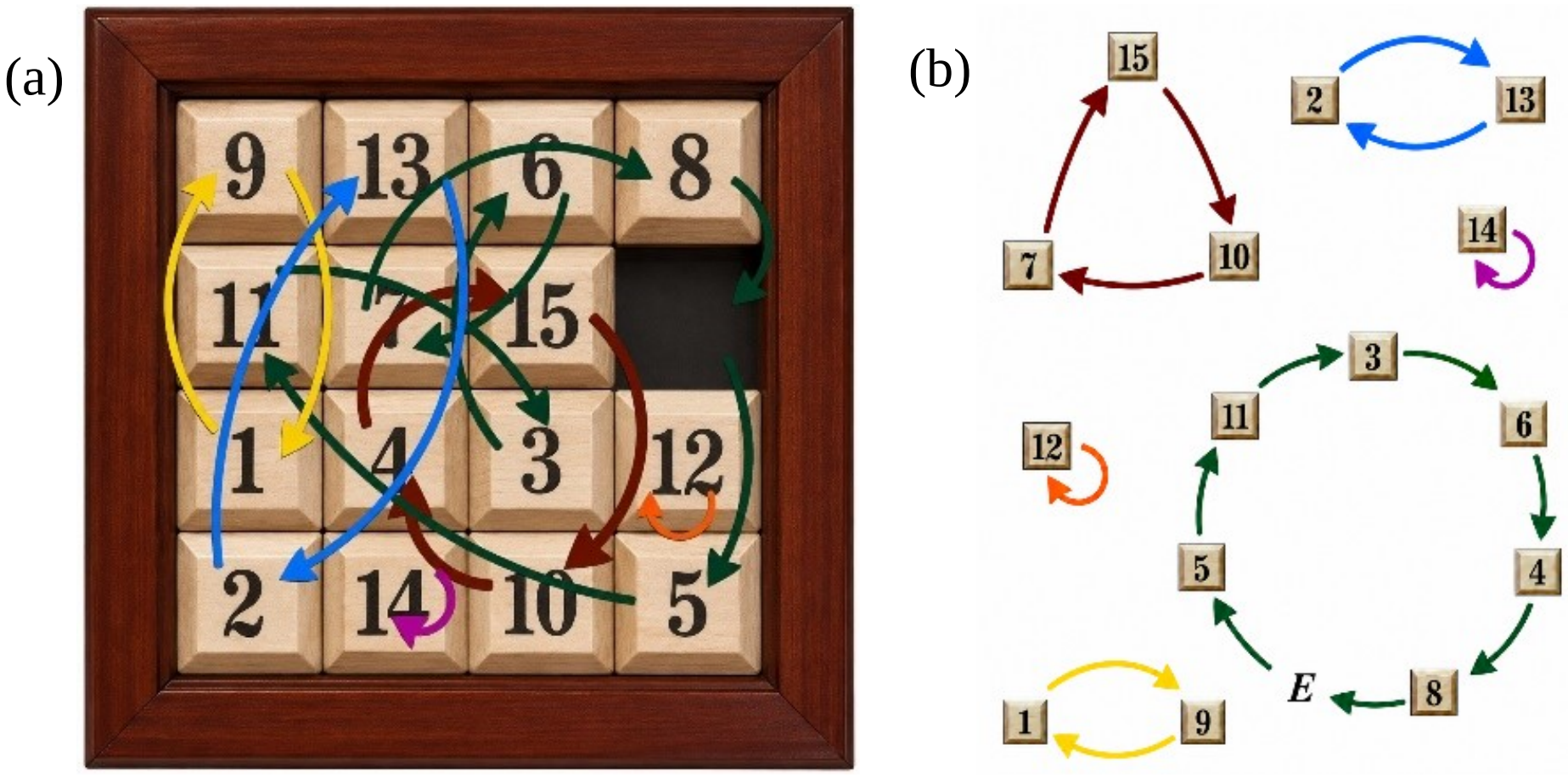

*Figure 2: From a puzzle configuration to its cycle representation: (a) the arrows shown on the puzzle; (b) the same arrows rearranged to show the cycles.*

A legal move swaps the moved tile with the empty element $E$. Consequently, only the cycle or cycles containing these two elements can change; all other cycles remain unchanged. Two cases are possible. If the moved tile and $E$ belong to different cycles, the move joins these cycles into a single cycle. If they belong to the same cycle, the move splits that cycle into two.

Figure 3 illustrates the first case. Before the move, the tile labeled 15 and the empty element $E$ belong to different cycles. Sliding tile 15 into the empty square changes only the two relevant connections in the arrow diagram (see the blue arrows in Figure 3), and the two original cycles merge into one. On the other hand, when the two elements belong to the same cycle, the legal move breaks this cycle into two separate cycles (see Figure 4).

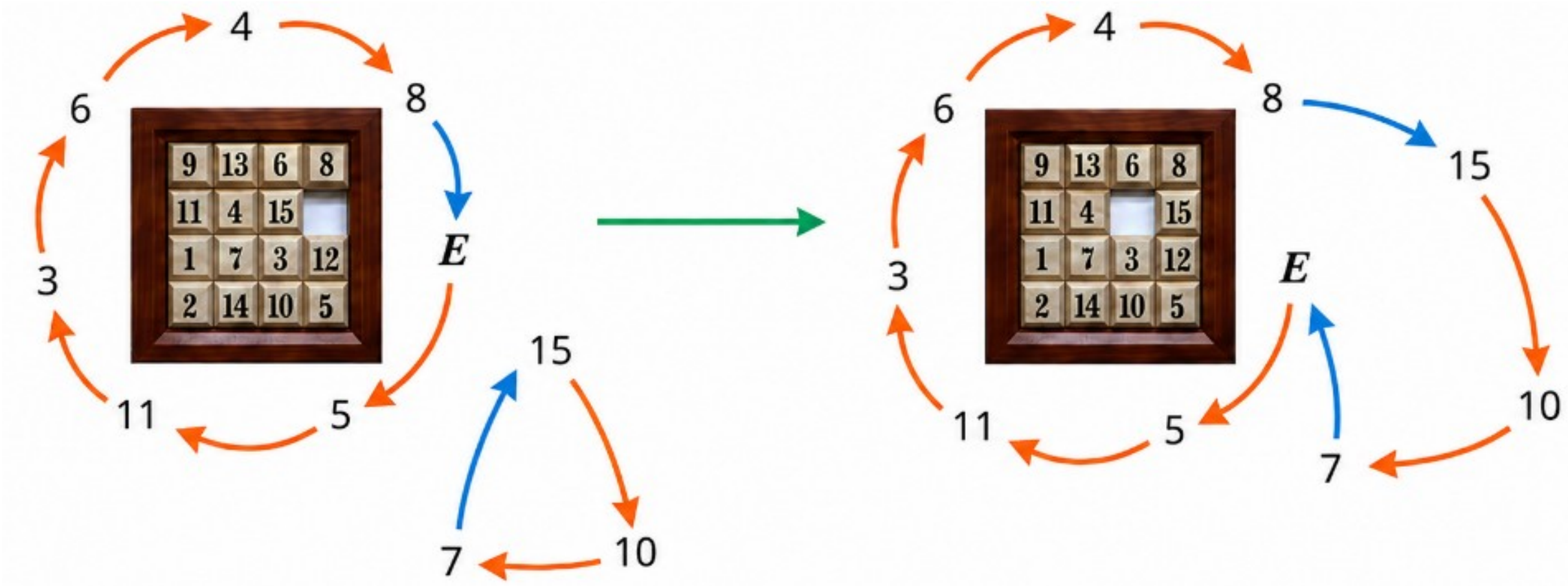

*Figure 3: The effect of a legal move on the cycle decomposition: two cycles merge.*

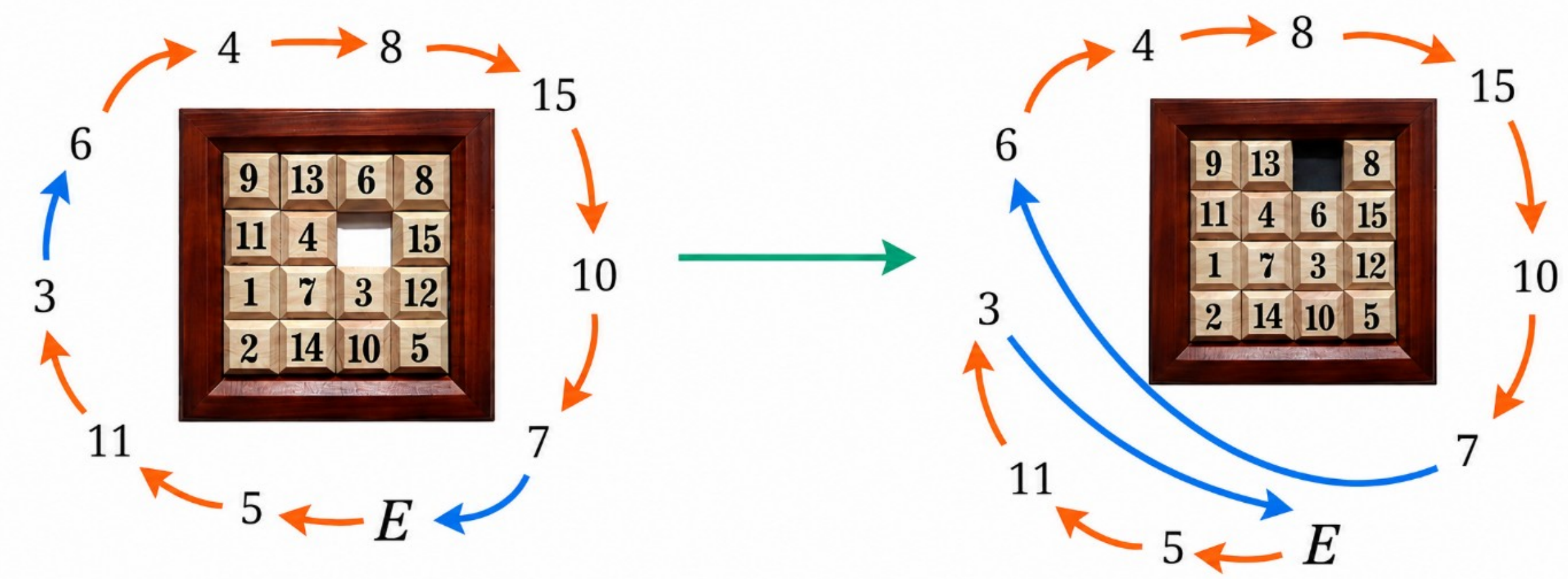

*Figure 4: The effect of a legal move on the cycle decomposition: one cycle splits into two.*

Every legal move therefore changes the number of cycles by exactly one: two cycles merge, or one cycle splits into two. In either case the parity of the total number of cycles changes.

**The Solvability Criterion**

Color the board in a checkerboard pattern so that adjacent squares have opposite colors and the bottom-right square is black. Since each legal move shifts a tile adjacent to the empty square, the square occupied by *E* changes color after every move. At the same time, each move changes the number of cycles by exactly one, and hence also changes its parity. The relationship between these two properties is therefore invariant. Consequently, any solvable configuration must have the same relationship between the parity of the number of cycles and the color of the empty square as the solved configuration.

In the solved configuration, *E* occupies the black bottom-right square. The solved configuration has sixteen one-element cycles and thus an even number of cycles. In the 15–14 configuration, the empty square remains on the same black square, but tiles 14 and 15 form a two-element cycle, giving fifteen cycles in total. The two configurations therefore have different parity-color relationships and cannot be connected by legal moves. Hence the 15–14 configuration is unsolvable.

With a little practice, it is relatively straightforward to place tiles 1 through 13 in their correct positions. This explains why early solvers could easily gain the impression that every starting configuration was solvable: most of the puzzle could indeed be completed before the obstruction became apparent. Once the first thirteen tiles are fixed, only tiles 14 and 15 and the empty square remain. Moving the empty square among the three remaining positions does not change the left-to-right order of tiles 14 and 15. If tile 14 is to the left of tile 15, the puzzle can be completed by sliding the two tiles into place. If tile 15 is to the left of tile 14, the configuration is connected by such moves to the unsolvable 15–14 configuration. Since legal moves preserve the invariant, an initial configuration satisfying the parity-color condition must lead to the first case and can be completed. Thus, the condition is sufficient as well as necessary.

Johnson and Story (1879) gave one of the earliest published mathematical treatments of the solvability criterion. Their paper is dated December 1879, the same year in which the puzzle was first marketed.

# METHOD

## Participants and Educational Context

The first author taught the 15-puzzle lesson as a visiting teacher in 16 classes at four schools in Szeged, Hungary. Table 1 gives the number of participants and classes at each grade level. The lesson aimed to teach students to construct a cycle representation of a puzzle configuration and use it to decide whether the configuration was solvable. The study examined how successfully students could do this immediately afterwards.

| **Grade** | ***n*** | **Number of classes** |
|---|---|---|
| 5 | 28 | 1 |
| 6 | 51 | 2 |
| 7 | 21 | 1 |
| 8 | 47 | 3 |
| 9 | 109 | 6 |
| 10 | 40 | 2 |
| 11 | 17 | 1 |

Table 1: Distribution of participants by grade level

Altogether, 313 students in Grades 5–11 attended the lessons. All of them submitted a response sheet, and all responses were included in the analysis. Each double-sided response sheet contained the questionnaire on one side and a puzzle task on the other. On the puzzle-task side, students were asked to construct a cycle representation of the given configuration and decide whether it was solvable. One version presented a solvable configuration and was given to 155 students; the other presented an unsolvable configuration and was given to 158. Within each class, the two versions were distributed as evenly as possible, with students sitting in alternating columns receiving different versions. The students did not know that two versions were in use.

### Lesson Design

All 16 lessons followed the same sequence. Students first had time to experiment with the puzzle, after which the unsolvable 15–14 configuration was presented. Arrow diagrams were then introduced as a hypothetical solution plan, showing where each tile should be placed if the tiles could be removed from the box, and students learned to record the resulting cycles in cycle notation. The final step was to use the parity of the number of cycles together with the color of the empty square to decide whether a configuration was solvable. Thus, the lesson moved from hands-on activity to visual representation and then to symbolic, structural reasoning. Each lesson lasted about 45 minutes; the approximate timing is shown in Table 2.

| **Approx. time** | **Phase** | **Main activity** | **Representation / resources** |
|---|---|---|---|
| 0–10 min | Introduction and exploration | The 15-puzzle was introduced, students were invited to try solving it, and brief strategic hints were provided when needed. | Physical and digital versions of the puzzle |
| 10–14 min | Historical background and problem posing | The history of the puzzle was briefly presented, followed by an example of a configuration that cannot be solved. | Slides and visual demonstration |
| 14–25 min | Mathematical analysis | Checkerboard coloring was introduced, followed by cycle representation. Students were shown how the number of cycles changes after a legal move and how this relates to the color of the square occupied by the empty element. The solvability criterion was then formulated. | Checkerboard coloring, arrow/cycle representations, classification rule |
| 25–32 min | Guided practice | A further configuration was analyzed together. Students constructed its cycles, identified the background color of the empty square, and used the classification rule to decide whether the configuration was solvable. | An example and the classification rule displayed to the class |
| 32–42 min | Individual assessment | Students received a double-sided response sheet. They completed the questionnaire side first and then solved the puzzle task on the other side using the instructions provided. | Double-sided response sheet |
| 42–45 min | Closing | Response sheets were collected, and students were informed that the two task versions contained different configurations, one solvable and one unsolvable. | Whole-class discussion |

Table 2: Main phases of the 15-puzzle lesson

Throughout the mathematical part of the lesson and the puzzle task, students could refer to the displayed classification table (see Figure A1 in the Appendix for the original Hungarian table and Table 3 for a simplified English version). They were not expected to derive or prove this rule on their own. The assessed task was to convert a concrete configuration into cycle notation and then use the supplied rule correctly.

| **Number of cycles** | **Background color of the empty square** | **Classification** |
| --- | --- | --- |
| Even | Black | Solvable |
| Even | White | Unsolvable |
| Odd | Black | Unsolvable |
| Odd | White | Solvable |

Table 3: Classification rule used during the lesson

**Data Collection**

At the end of the lesson, students completed the questionnaire side of the response sheet and then, when instructed, turned to the puzzle task on the other side. The questionnaire was in Hungarian (see Figure A2 in the Appendix); the English translations of its four five-point Likert-type items are listed below. Response options were 1 (not at all), 2 (rather not), 3 (mixed), 4 (rather yes), and 5 (completely).

- Did you find the lesson interesting?
- Were you surprised that mathematics also has such playful aspects?
- Would you have thought that such a simple game has such a substantial mathematical background?
- How well did you understand what was presented?

An optional space for comments was also provided, but these comments were not part of the quantitative analysis. The four rating items were intended to capture different reactions to the lesson, so they were analyzed separately rather than treated as a scale. In particular, high ratings on the questions about surprise and prior expectations were not interpreted as stronger approval; those items were used only to describe students' immediate reactions.

The second part asked students to analyze a given 15-puzzle configuration (see Figure A3 in the Appendix). Two versions of the task were used: one contained a solvable configuration and the other an unsolvable configuration. The two configurations were deliberately chosen to be very similar (Figure 5). The empty square occupied the same position in both configurations, and most of the cycles were identical. The configurations differed only in how five elements were grouped into cycles, resulting in opposite parities of the total number of cycles and hence opposite solvability classifications. In both versions, students were required first to construct the cycle representation of the configuration and then to decide whether it was solvable. The solved configuration was shown next to the configuration to be analyzed as a reference. The worksheet also included a short reminder of the intended sequence: construct the cycles, check that all elements are represented, determine the number of cycles, examine the background color of the empty square, and decide whether the configuration is solvable.

(a)

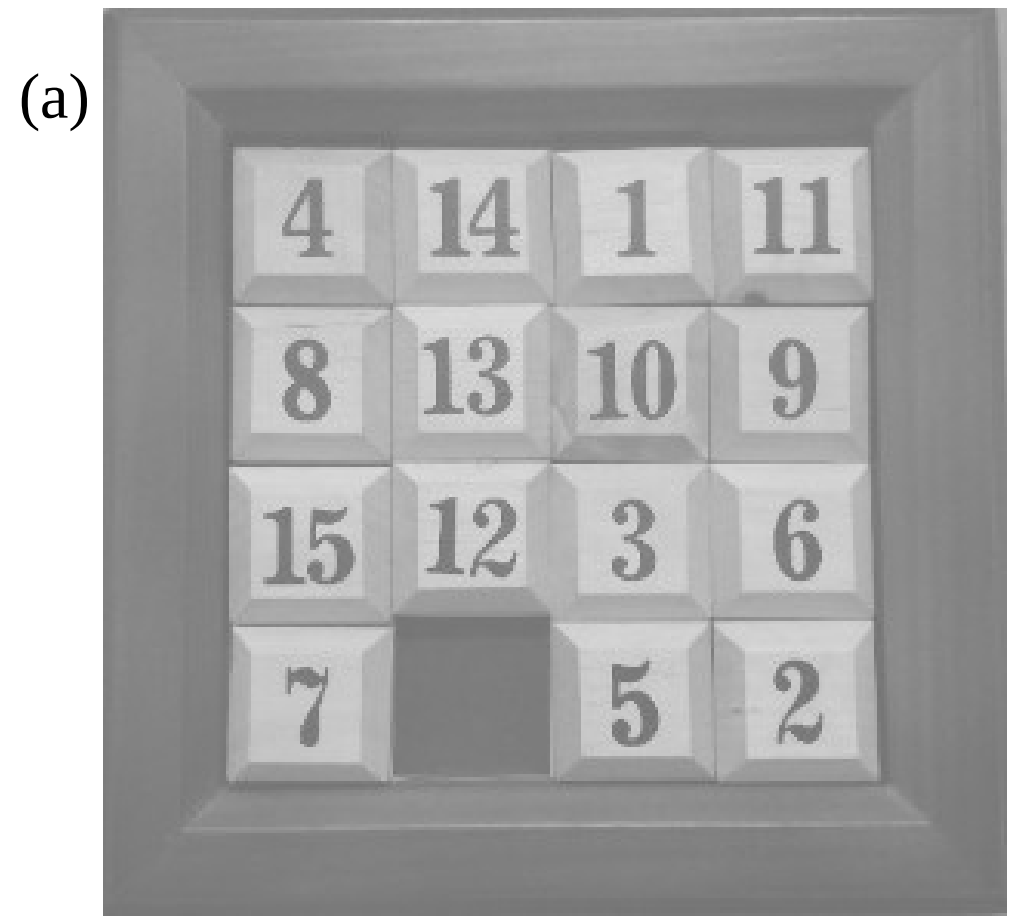

(b)

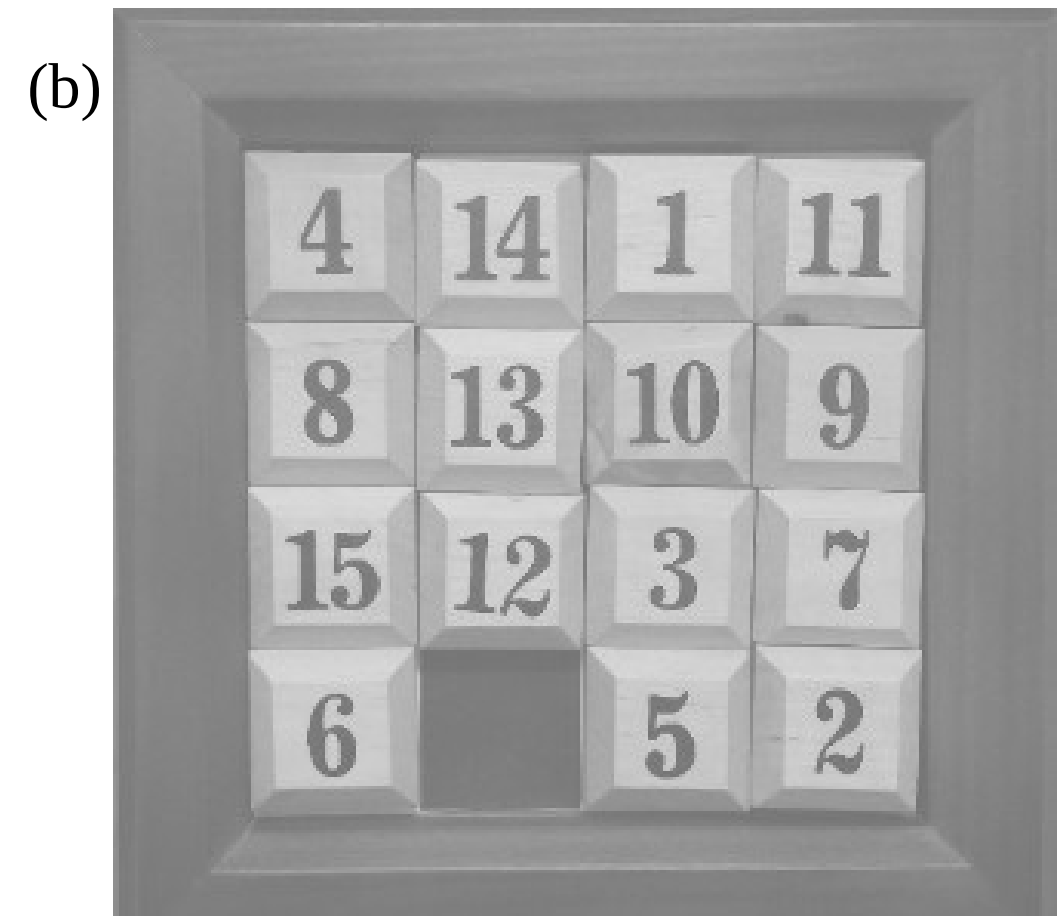

Figure 5: The two puzzle configurations used in the worksheet: (a) solvable configuration; (b) unsolvable configuration.

*Note.* The correct cycle representation for configuration (a) was
(1 4 11 3) (2 14 *E*) (5 8 9 15) (6 13 7 10 12), and the configuration was classified as solvable.
The correct cycle representation for configuration (b) was
(1 4 11 3) (2 14 *E*) (5 8 9 15) (6 13) (7 10 12), and the configuration was classified as unsolvable.

The classification table (see Table 3) remained visible while students completed the puzzle task. Accordingly, students were not asked for a written proof or a full justification of the solvability criterion. The task focused on a more immediate question: could they construct the cycles for a new configuration and use the displayed rule to decide whether it was solvable?

### Coding and Statistical Analysis

Each Likert-type item of the questionnaire was summarized separately by its median, interquartile range, and the percentage of students choosing 4 or 5. Because the items served different descriptive purposes, no total score was formed.

Worksheet responses were coded by the first author. The cycle representation was counted as correct only when all 16 elements, including the empty square, appeared exactly once and the decomposition matched the given configuration. Equivalent forms obtained by changing the order of disjoint cycles or choosing a different starting point within a cycle were accepted. For the final classification, a response was coded as correct only when the student recorded a solvability decision and that decision was correct; an omitted final classification was therefore coded as not correct. A combined outcome recorded whether both the cycle representation and the final classification were correct.

As a supplementary descriptive measure, we also recorded whether the final classification was consistent with the student's own written work. This was assessed whenever the number of cycles intended by the student could be identified from the cycle or cycle-like representation on the worksheet, even if that representation was mathematically incorrect and therefore received no credit for cycle construction. In five cases, students incorrectly labeled the parity of their cycle count or the color of the empty square but then used the value they had recorded consistently. Since this supplementary coding was intended to capture the subsequent use of the recorded information rather than knowledge of parity or color identification, these responses were treated as consistent. Students who gave no final classification were not scored on this supplementary variable.

We report the number and percentage of students with a correct cycle representation, a correct solvability classification, and both components correct. Students who did not record a final classification were included in the denominator and were counted as not reaching a correct classification. The supplementary consistency measure was reported descriptively among students who provided a final classification.

Differences between the solvable and unsolvable task versions were examined using Cochran–Mantel–Haenszel tests stratified by class, reflecting the within-class allocation of the two worksheet versions. Separate comparisons were made for correct cycle representation, correct classification, and the combined outcome. Mantel–Haenszel common odds ratios were reported as summary effect estimates, together with 95% confidence intervals. The three corresponding $p$-values were adjusted together using Holm's procedure. Classification accuracy among students who had constructed the cycle representation correctly was examined descriptively only because this subgroup was defined by post-assignment performance.

The examples of incorrect representations presented in the Appendix were selected illustratively after the quantitative coding and were not treated as additional outcome variables.

All tests were two-sided. Analyses were carried out in R version 4.5.1.

## Ethical Considerations

The lessons formed part of regular classroom instruction and were conducted with the participating mathematics teachers' permission and in consultation with them. Worksheets were anonymous: students did not provide names or other information that could be used to link responses to individual students. Participation in the activity had no effect on grades and did not place students at an educational disadvantage. Neither students nor their teachers later received class-level performance results. The lesson also did not replace or reduce students' ordinary opportunities to learn mathematics.

# RESULTS

## Student Responses to the Lesson

Responses to the four questionnaire items are summarized in Table 4. Interest in the lesson and self-reported understanding received high ratings, whereas the two items concerning surprise and prior expectations produced more middle-range responses.

| Questionnaire item | Median $[Q_1, Q_3]$ | Responses 4–5 (%) |
|---|---|---|
| Interest in the lesson | 5 [4, 5] | 91.1 |
| Surprise at the playful aspects of mathematics | 3 [2, 4] | 37.1 |
| Prior expectation of substantial mathematical background | 3 [2, 4] | 31.3 |
| Self-reported understanding | 5 [4, 5] | 83.7 |

Table 4: Descriptive statistics for the questionnaire items

## Worksheet Performance

Of the 313 students, 245 (78.3%) produced the correct cycle representation. A correct final solvability classification was recorded for 225 students (71.9%), and 212 (67.7%) were correct on both parts. A final classification was not recorded by 41 students (13.1%).

The incorrect cycle representations took several forms. Examples included incomplete or cycle-like constructions that did not form a valid decomposition, omission of one or more elements, duplication of a cycle or some of its elements, and an incorrect connection in an otherwise largely correct representation. Illustrative examples are shown in the Appendix (Figures A4–A7).

Of the 245 students with a correct cycle representation, 212 (86.5%) also reached the correct final classification. Among the 68 students whose cycle representation was incorrect, 13 (19.1%) nevertheless gave the correct final classification, while 33 did not record a final classification.

The supplementary coding gives a more detailed view of how students used their own written work. Among the 272 students who recorded a final classification, 241 (88.6%) gave a classification consistent with applying the displayed rule to the values recorded on their worksheet. This was the case for 218 of 237 students (92.0%) with a correct cycle representation and 23 of 35 (65.7%) with an incorrect one. Thus, students with an incorrect cycle representation were also less likely to apply the classification rule consistently. Because only 35 students were in the latter group and the grouping was based on students' performance, this comparison is descriptive and does not explain why the difference occurred.

Table 5 gives the results by grade for descriptive purposes. The percentage of correct cycle representations ranged from 50.0% in Grade 5 to 84.4% in Grade 9, while the percentage with both parts correct ranged from 46.4% to 75.2%. Several upper grades showed higher cycle-construction accuracy, but the pattern was not monotonic. The performance across the age range is particularly noteworthy because cycle notation is not part of the Hungarian school curriculum: 73.5% of students in Grades 5–8 and 82.5% of those in Grades 9–11 produced a correct cycle representation after the lesson. Since grade and class membership were partly confounded in this sample, these figures should not be read as estimates of a grade-level effect.

| **Grade** | *n* | **Correct cycle representation, *n* (%)** | **Correct representation and classification, *n* (%)** |
|---|---|---|---|
| 5 | 28 | 14 (50.0%) | 13 (46.4%) |
| 6 | 51 | 40 (78.4%) | 34 (66.7%) |
| 7 | 21 | 15 (71.4%) | 13 (61.9%) |
| 8 | 47 | 39 (83.0%) | 32 (68.1%) |
| 9 | 109 | 92 (84.4%) | 82 (75.2%) |
| 10 | 40 | 31 (77.5%) | 26 (65.0%) |
| 11 | 17 | 14 (82.4%) | 12 (70.6%) |
| Total | 313 | 245 (78.3%) | 212 (67.7%) |

Table 5: Worksheet performance by grade level

### Differences Between the Two Task Versions

Table 6 compares performance on the solvable and unsolvable task versions using analyses stratified by class. Cycle-representation accuracy was similar for the two versions. Correct classification was more frequent for the solvable configuration, and this difference remained statistically significant after Holm correction. The combined outcome was also more frequent for the solvable version, but this difference did not remain statistically significant after correction.

The final classification was left blank by 15 students (9.7%) on the solvable version and 26 (16.5%) on the unsolvable version.

Among students who had first constructed the cycles correctly, 113 of 124 (91.1%) on the solvable version and 99 of 121 (81.8%) on the unsolvable version also reached the correct classification. We treated this as a descriptive comparison only, because entry into this subgroup depended on students' post-assignment performance.

| **Outcome** | **Solvable task** | **Unsolvable task** | $OR_{MH}$ **[95% CI]** | $p$ | $p_{Holm}$ |
|---|---|---|---|---|---|
| Correct cycle representation | 124/155 (80.0%) | 121/158 (76.6%) | 1.258 [0.724, 2.186] | 0.416 | 0.416 |
| Correct classification | 122/155 (78.7%) | 103/158 (65.2%) | 1.980 [1.195, 3.280] | 0.007 | 0.022 |
| Correct representation and classification | 113/155 (72.9%) | 99/158 (62.7%) | 1.625 [1.003, 2.632] | 0.049 | 0.098 |

Table 6: Class-stratified comparisons of the solvable and unsolvable task versions

*Note.* Cochran–Mantel–Haenszel analyses were stratified by class. Odds ratios greater than 1 favor the solvable task version. The three *p*-values were adjusted together using Holm's procedure. The reported 95% confidence intervals are conventional, unadjusted confidence intervals.

## DISCUSSION

The clearest result is that most students were able to produce the cycle representation after one 45-minute, teacher-guided lesson. Of 313 participants, 245 (78.3%) did so correctly, and 212 of these 245 students (86.5%) also reached the correct classification. For an introductory lesson spanning Grades 5–11, this suggests that cycle notation was accessible to a substantial majority of the participating students. Correct cycle representations were more common in several of the upper grades, although the pattern was irregular. Because grade and class could not be disentangled in the present design, we do not interpret this as evidence of a developmental trend.

The use of several representations in the lesson can be viewed in light of Dienes's principle of perceptual variability, research on multiple representations, and Sfard's distinction between operational and structural views of mathematical objects (Ainsworth, 2006; Dienes, 1960; Duval, 2006; Sfard, 1991). The board, arrow diagram, and cycle notation represented the same puzzle configuration in different forms. Cycle notation condensed the configuration into a form in which a global property relevant to solvability could be examined. From a mathematical point of view, these representations placed a finite, discrete function in a concrete setting, although neither functions nor mappings were discussed with the students. Our data therefore do not show that the lesson broadened students' concept of function. They do, however, show that many students could work with the cycle representation itself, which makes that longer-term question worth investigating.

The results also show why cycle construction and classification should not be treated as the same outcome. Cycle-construction accuracy did not differ clearly between the two task versions, whereas correct final classification was more frequent for the solvable configuration. Even among students who had constructed the cycles correctly, classification was less often correct for the unsolvable task. The difficulty, then, does not appear to lie simply in producing the notation. It may arise when students have to connect the parity of the number of cycles, the color of the empty square, and the final entry in the classification rule. The present data do not tell us which of these steps was responsible for the difference.

One possible explanation is that the two conclusions are not psychologically symmetric. The higher number of omitted final classifications on the unsolvable version may also have contributed to the difference in the primary classification outcome. Declaring a configuration solvable may be easier than accepting that an arrangement is impossible. This remains only a possibility, because students were not asked to explain their reasoning. Written justifications or interviews would be needed to determine whether errors came from counting cycles, coordinating the parity of the number of cycles with the color of the empty square, reading the classification table, or from the negative conclusion itself.

The supplementary coding adds some detail to the distinction between cycle construction and classification. Among the 35 students with an incorrect cycle representation who nevertheless gave a final classification, 23 gave a classification consistent with applying the displayed rule to the values recorded on their worksheet. Since the classification table remained visible throughout the task, however, this should not be taken as evidence that these students had independently mastered the solvability criterion.

For teaching, this points to the importance of discussing what the representation is used for, not only how it is written. Placing solvable and unsolvable configurations side by side, asking students to state the criterion in their own words, and requiring a short justification could help make the decisive structural difference more visible. Such comparisons can draw attention away from surface appearance and towards the invariant relation between the parity of the number of cycles and the background color of the empty square.

## LIMITATIONS

The study has several limitations. There was no control group, pre-test, or alternative instructional condition. The results therefore describe what students were able to do immediately after this particular lesson; they do not show how much learning occurred, whether the puzzle-based approach outperformed another method, or how much of the learning would be retained.

The sample consisted of existing classes rather than students drawn randomly from a wider population, and students were clustered within classes. We accounted for class in comparisons between the two worksheet versions, but grade-level results were kept descriptive because grade and class membership were partly confounded. They should not be interpreted as general age or grade effects. The first author also taught all of the lessons and carried out the study. Replication with other teachers and in other schools would therefore be valuable.

The questionnaire is another limitation. In particular, the item on understanding measured students' own perception rather than conceptual knowledge. A high rating may have reflected confidence or enjoyment as well as understanding. Interviews, explanation tasks, and delayed assessment would provide stronger evidence about what students understood and whether it was retained.

The worksheet also assessed a deliberately limited form of performance. Students were asked to apply a representation and a classification rule that had just been taught and remained visible. They did not have to derive the criterion or justify it in writing. A correct answer should therefore not be taken as evidence of complete structural understanding, an independently constructed proof, or transfer to an unfamiliar context. The supplementary consistency coding has a similar limitation: it records whether a final classification fits the intermediate values written on the worksheet, not whether the student consciously derived the answer from them.

Finally, each version contained only one configuration. Any difference between the solvable and unsolvable tasks may therefore be partly specific to those particular items. A stronger design would use several matched configurations and collect richer information about intermediate reasoning and typical errors, for example through written explanations or interviews.

## CONCLUSION

After one teacher-guided lesson, most of the participating students in Grades 5–11 were able to represent a 15-puzzle configuration using permutation cycles. This result suggests that cycle notation can be introduced meaningfully before advanced or university-level mathematics when it is supported by a concrete puzzle and intermediate visual representations.

At the same time, producing the cycles and using them were not the same achievement. Cycle-construction accuracy did not differ clearly between the two task versions, whereas the unsolvable configuration was classified correctly less often. Instruction should therefore make students explain what the cycles tell them and how the representation leads to a solvable or unsolvable conclusion, rather than stopping once the decomposition has been written.

A useful next step would be to repeat the study with several matched configurations, delayed assessment, and qualitative data on students' reasoning. Over a longer timescale, it would also be valuable to examine whether sustained work with permutation cycles helps students view functions more broadly by adding finite, discrete, and non-formulaic mappings to the examples they know.

## ACKNOWLEDGMENTS

The first author was supported by the Ferenc Móricz Doctoral Foundation. The second author was supported by the National Research, Development and Innovation Office of Hungary under grant ADVANCED 153383.

ChatGPT (OpenAI) was used during the preparation of this manuscript to suggest wording, reformulate some passages, and comment on the English, clarity, organization and consistency of the text. The authors critically reviewed and revised all AI-assisted text and made the final decisions about the mathematical and educational content. The authors take full responsibility for the final manuscript.

## APPENDIX

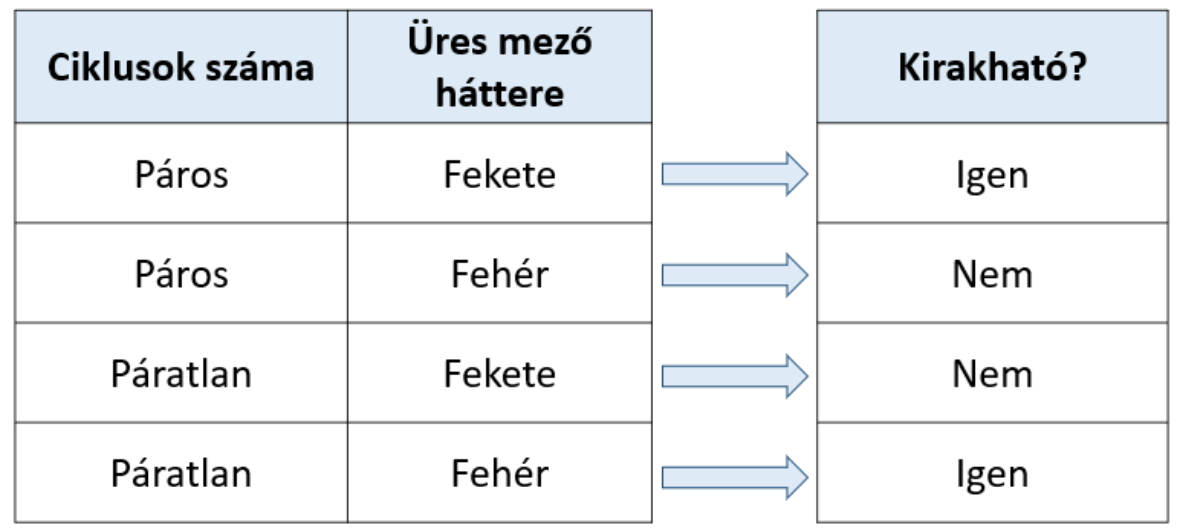

| Ciklusok száma | Üres mező háttere | | Kirakható? |
|---|---|---|---|
| Páros | Fekete | ⇨ | Igen |
| Páros | Fehér | ⇨ | Nem |
| Páratlan | Fekete | ⇨ | Nem |
| Páratlan | Fehér | ⇨ | Igen |

Figure A1: Original Hungarian classification table displayed during the lesson and kept visible during the worksheet task.

Osztály: ……………………………

***15-ös játék***

***Kérlek, jelezd minden kérdésnél a megfelelő szám bekarikázásával, hogy milyen mértékben jellemzőek Rád ezek az állítások.***

*1: egyáltalán nem*
*2: inkább nem*
*3: is-is*
*4: inkább igen*
*5: teljes mértékben*

| | | | | | |
|---|---|---|---|---|---|
| Érdekesnek találtad az órát? | 1 | 2 | 3 | 4 | 5 |
| Meglepett téged, hogy a matematikának vannak ilyen játékos részei is? | 1 | 2 | 3 | 4 | 5 |
| Gondoltad volna, hogy egy ilyen egyszerű játéknak ilyen komoly matematikai háttere van? | 1 | 2 | 3 | 4 | 5 |
| Mennyire értetted az elhangzottakat? | 1 | 2 | 3 | 4 | 5 |

*Egyéb észrevétel:*

***Kérlek, ne fordítsd meg a lapot!***

Figure A2: Original Hungarian post-lesson questionnaire completed before the individual worksheet task.

***Ki lehet-e rakni az első képen szereplő állásból a játékot?***
***Segítségképpen a második képen a kirakott játék szerepel.***
(Emlékeztető: ciklusok felrajzolása, szerepel-e az összes elem a rajzban, ciklusok számának meghatározása, háttérszín vizsgálata, kirakható-e)

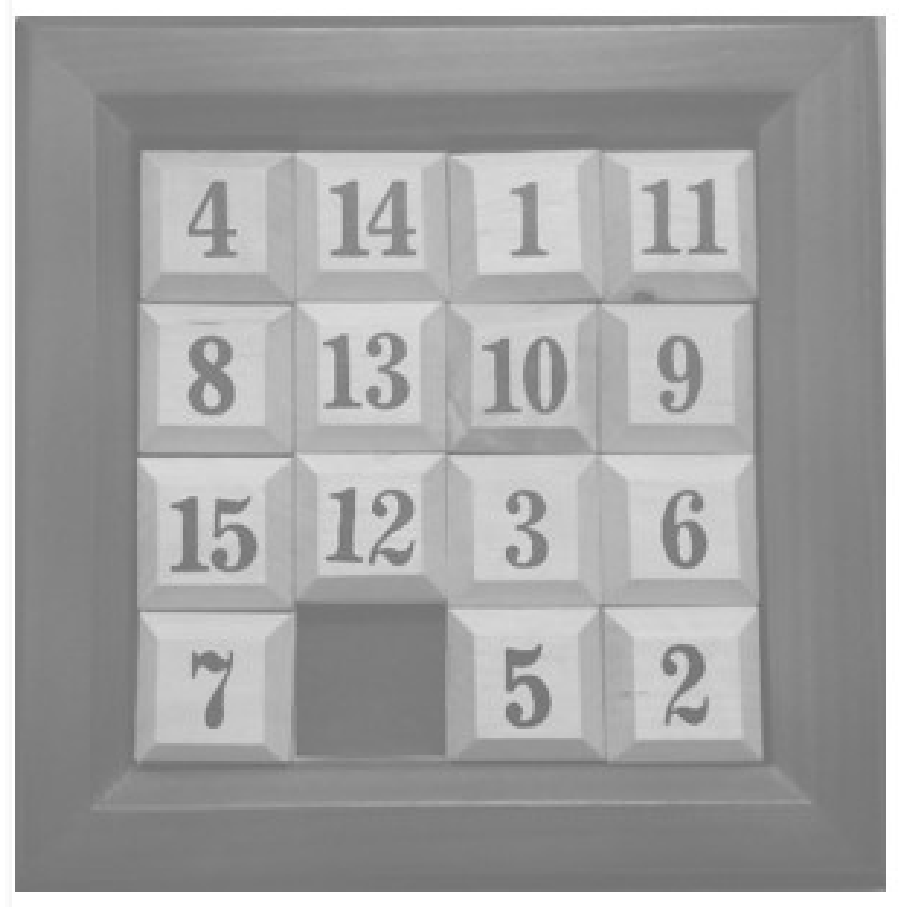

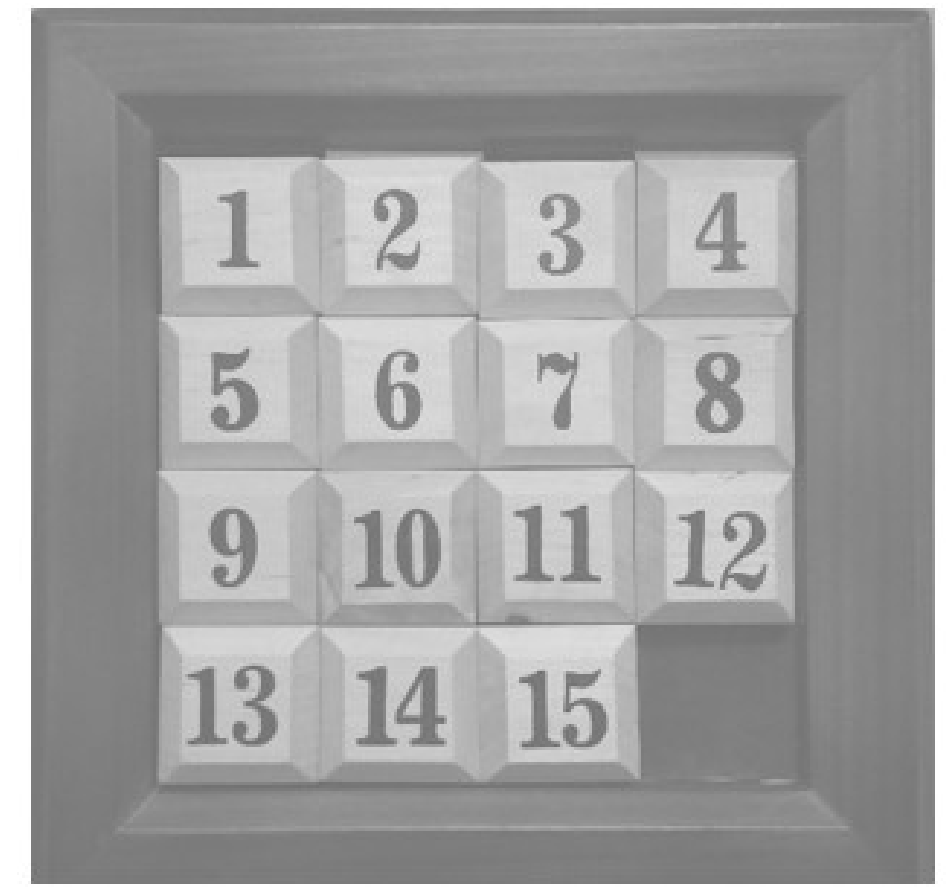

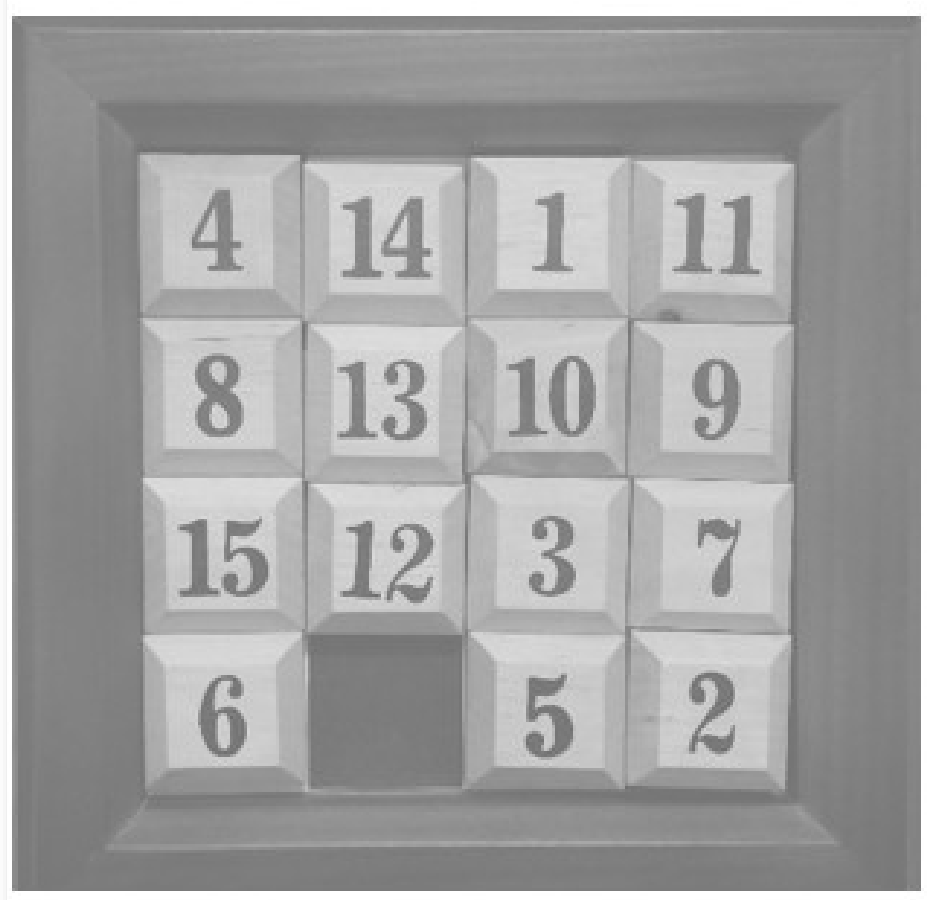

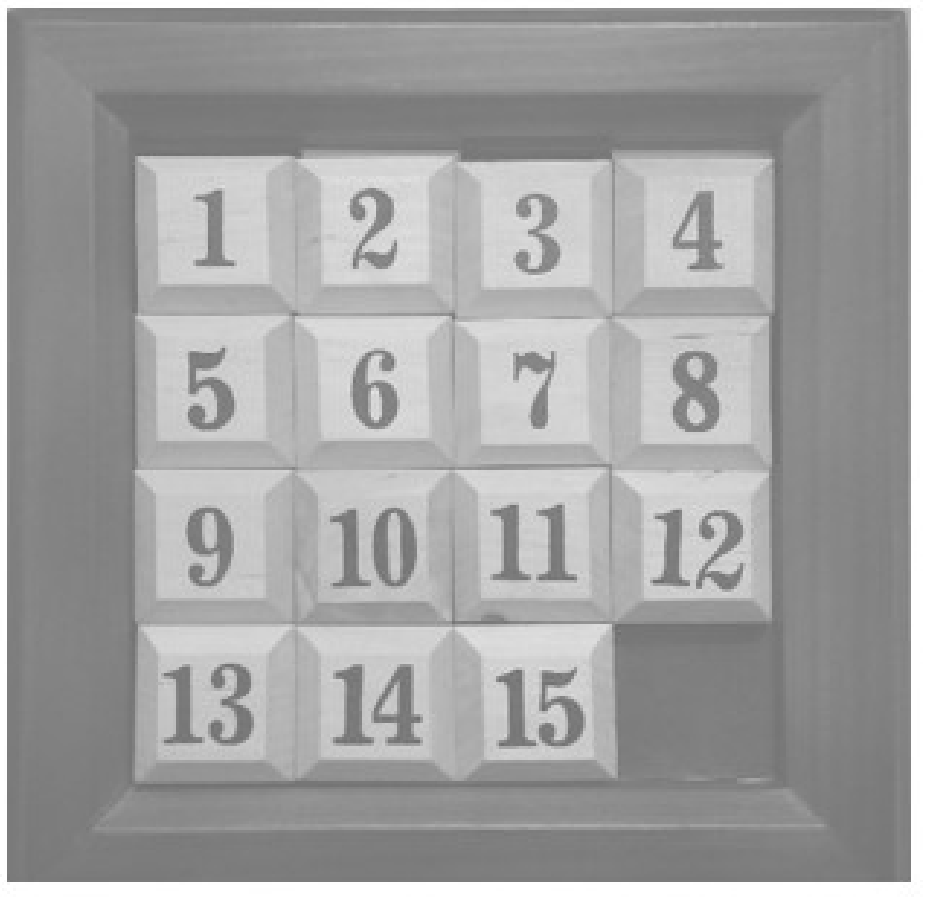

Figure A3: Original Hungarian worksheet versions used in the study: the solvable task version (above) and the unsolvable task version (below). In both versions, the solved configuration was shown alongside the configuration to be analyzed as a reference.

***Ki lehet-e rakni az első képen szereplő állásból a játékot?***
***Segítségképpen a második képen a kirakott játék szerepel.***
(Emlékeztető: ciklusok felrajzolása, szerepel-e az összes elem a rajzban, ciklusok számának meghatározása, háttérszín vizsgálata, kirakható-e)

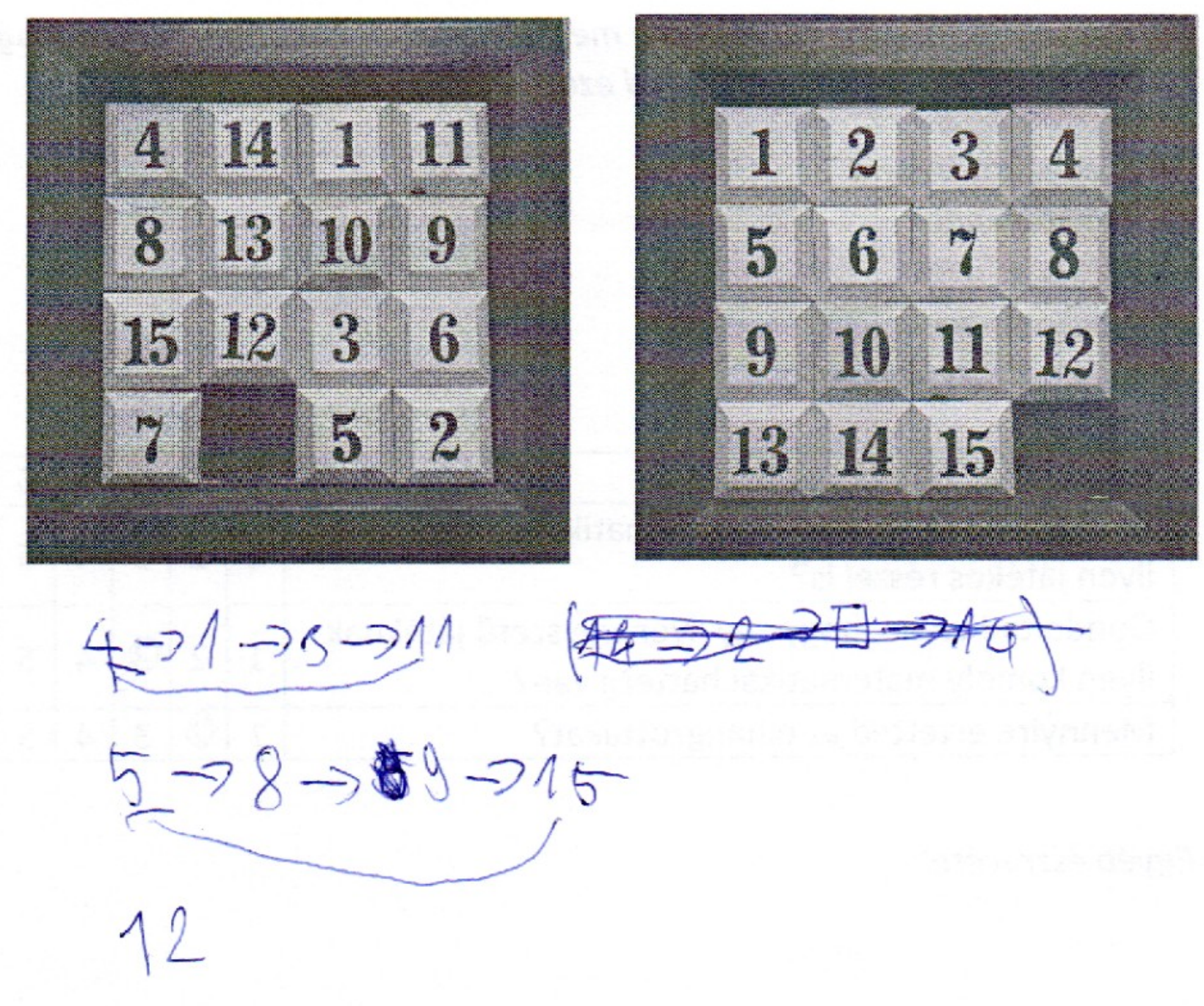

Figure A4: Anonymized student response illustrating a cycle representation in which one cycle was constructed correctly, while another was written in the inverse direction.

***Ki lehet-e rakni az első képen szereplő állásból a játékot?***
***Segítségképpen a második képen a kirakott játék szerepel.***
(Emlékeztető: ciklusok felrajzolása, szerepel-e az összes elem a rajzban, ciklusok számának meghatározása, háttérszín vizsgálata, kirakható-e)

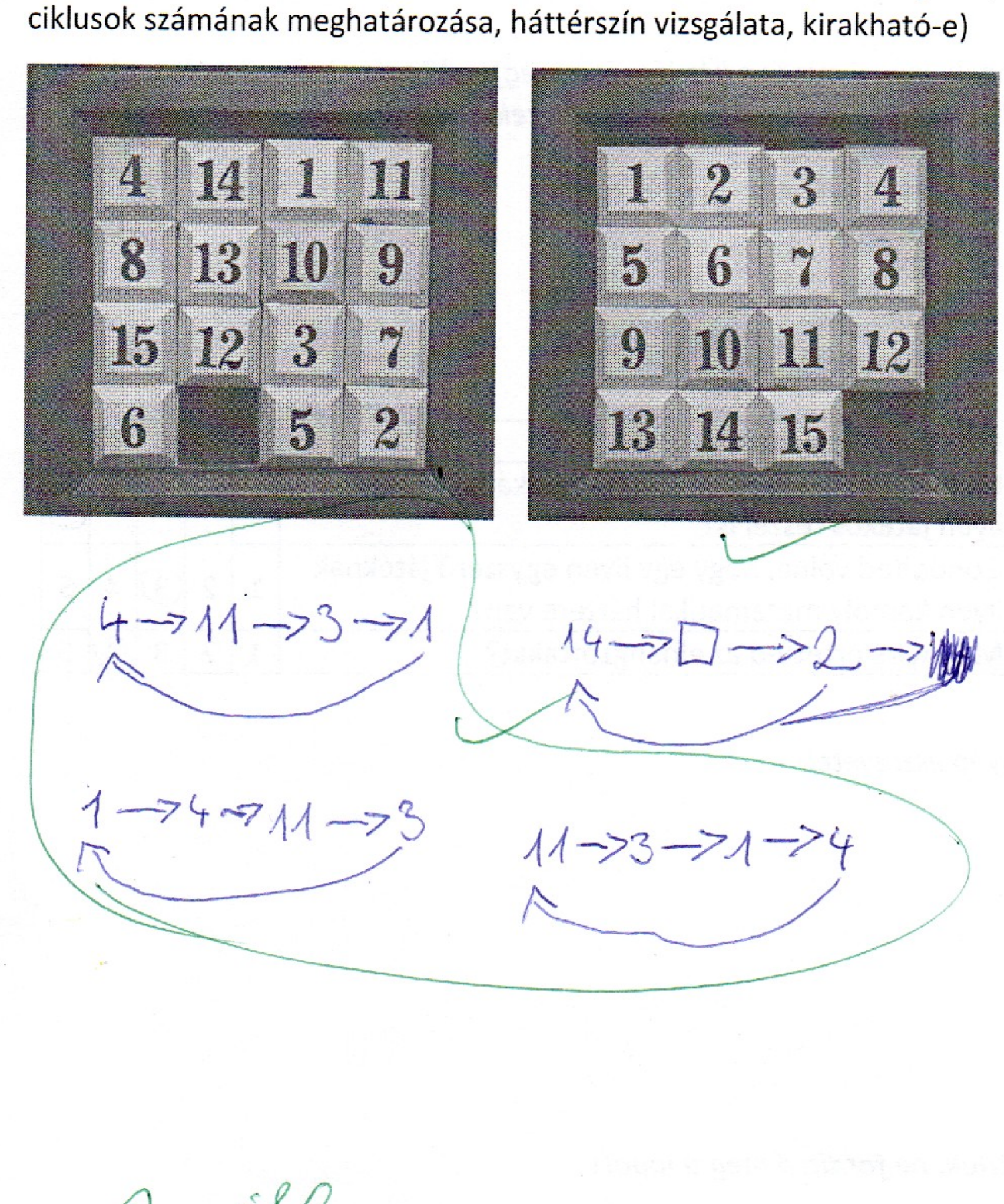

Figure A5: Anonymized student response illustrating duplication of the same cycle using different starting points.

***Ki lehet-e rakni az első képen szereplő állásból a játékot?***
***Segítségképpen a második képen a kirakott játék szerepel.***
(Emlékeztető: ciklusok felrajzolása, szerepel-e az összes elem a rajzban, ciklusok számának meghatározása, háttérszín vizsgálata, kirakható-e)

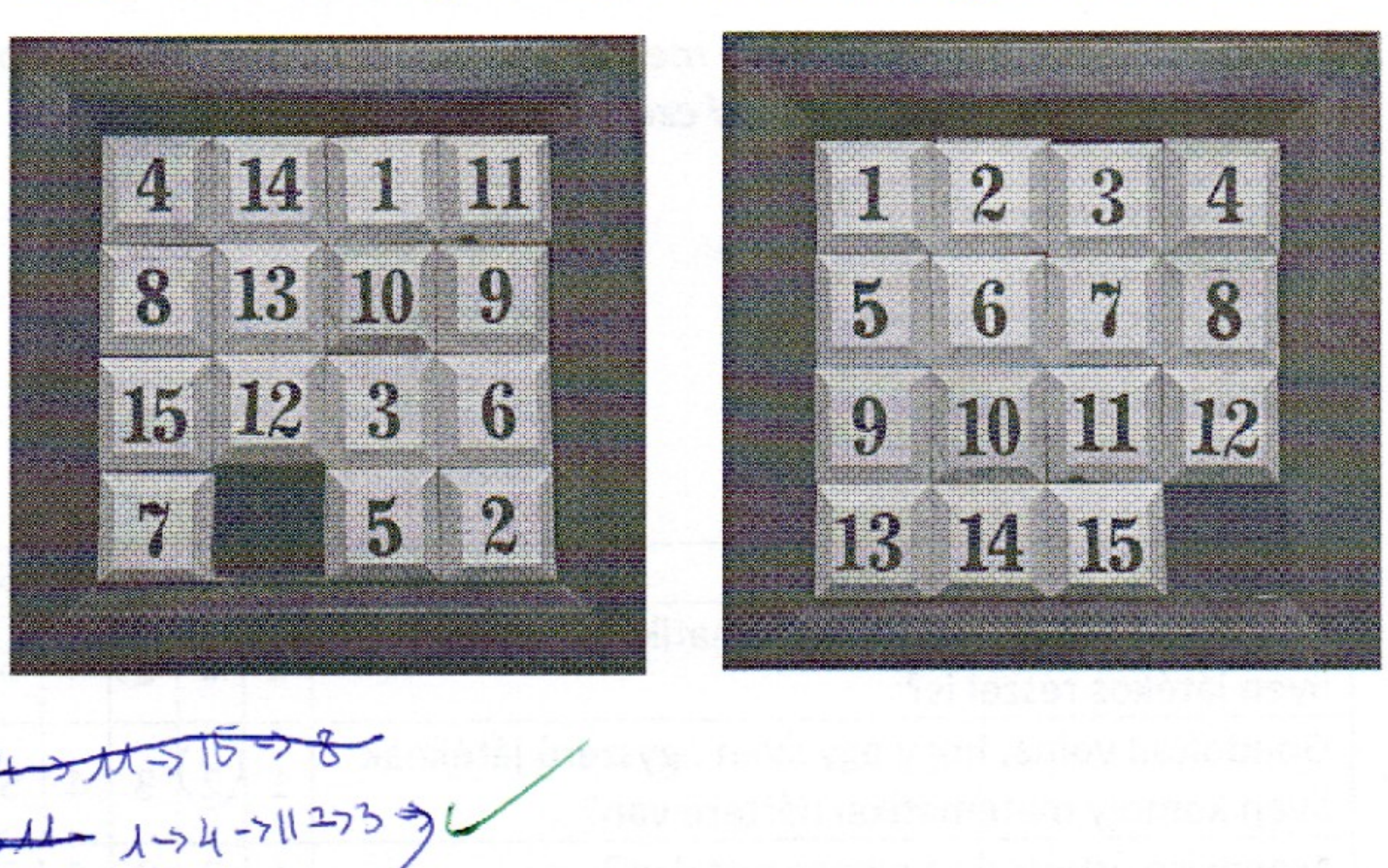

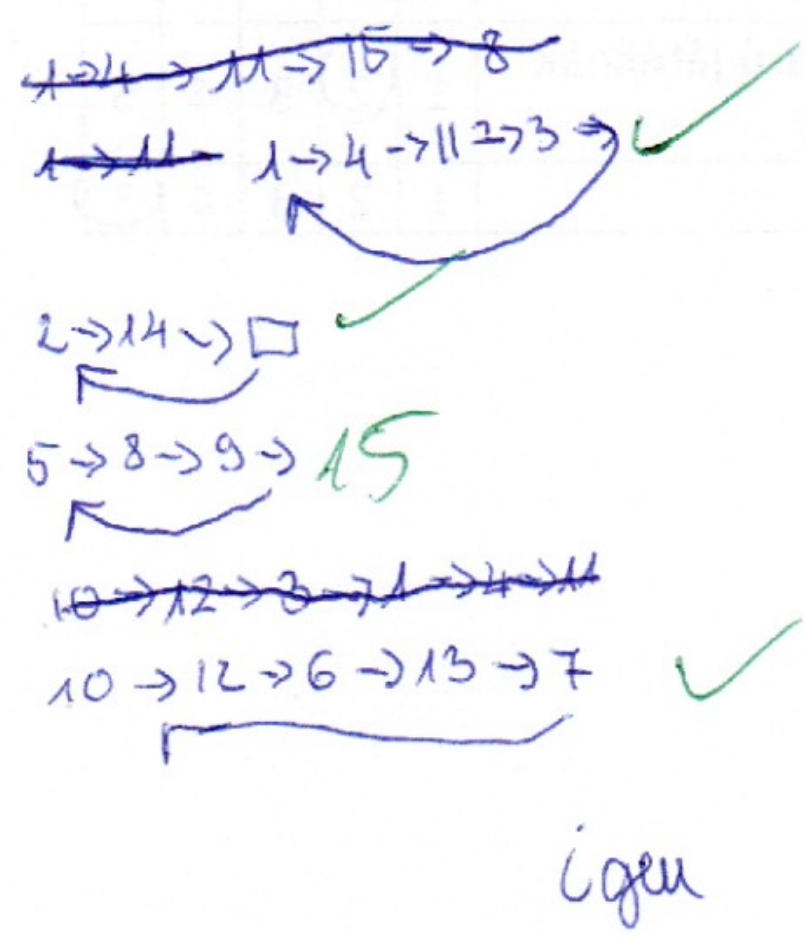

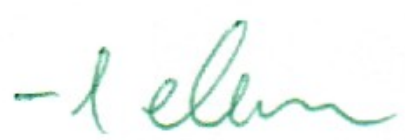

Figure A6: Anonymized student response illustrating omission of an element from an otherwise largely correct cycle construction.

***Ki lehet-e rakni az első képen szereplő állásból a játékot?***
***Segítségképpen a második képen a kirakott játék szerepel.***
(Emlékeztető: ciklusok felrajzolása, szerepel-e az összes elem a rajzban, ciklusok számának meghatározása, háttérszín vizsgálata, kirakható-e)

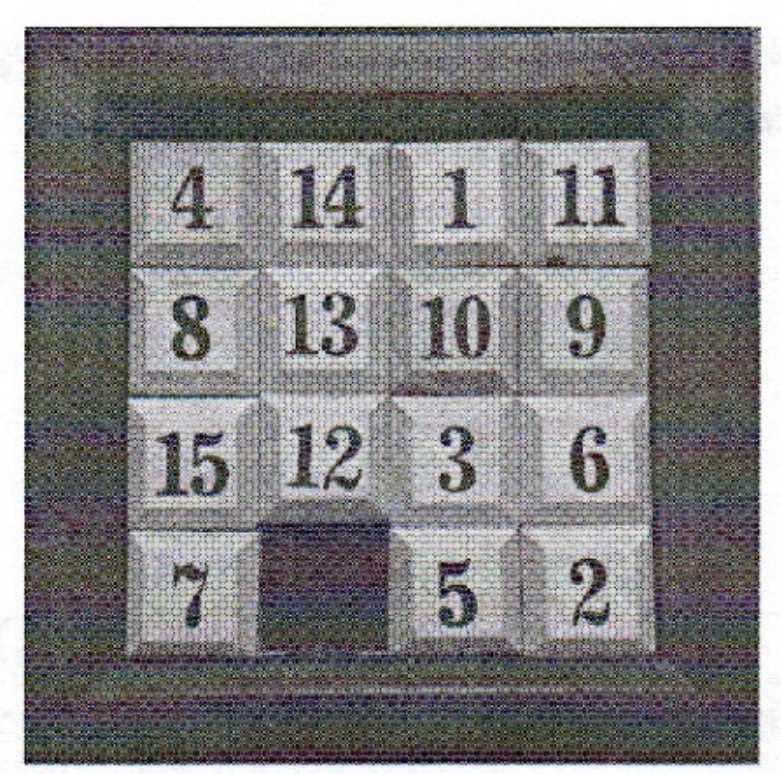

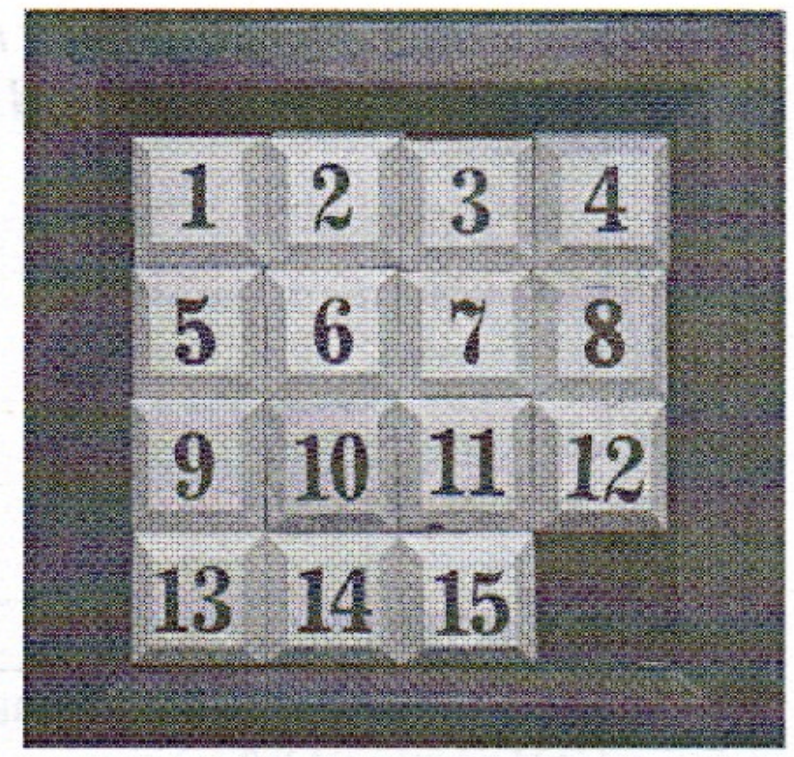

(1,4)(2,14)(3,1)(5,8,15)(6,13)

4→11 1→4 2→14 3→1
5→8 6→13 7→10 8→9
10→15 11→12 13←6 14→7
15→2

?

Figure A7: Anonymized student response in which several element-to-element mappings were recorded, with some errors, but the mappings were not successfully combined into a cycle decomposition.